\documentclass[11pt]{article}

	\usepackage[pdftex]{graphicx}
 
\usepackage{amssymb,amsmath}
\usepackage{hyperref}
\usepackage{url}

\def\R{\mathbf{R}}
\def\Z{\mathbf{Z}}
\def\H{\mathbf{H}}
\def\PP{P}

\def\proof{\par\medskip\noindent {\bf Proof.}{\hskip1em}}
\def\qed{~\vrule height .9ex width .8ex depth -.1ex}

\def\eps{\epsilon}
\def\vol{\mathrm{vol}}
\def\lng{\mathrm{length}}
\def\area{\mathrm{area}}

\newtheorem{thm}{Theorem}
\newtheorem{lem}{Lemma}
\newtheorem{cor}[thm]{Corollary}
\newtheorem{prop}[lem]{Proposition}
\newtheorem{defi}[lem]{Definition}

\newtheorem{ques}{Question}

\title{A discontinuous isoperimetric profile for a complete Riemannian manifold}
\author{Stefano Nardulli and Pierre Pansu}
\begin{document}

\maketitle

\begin{quote}
{\small 
ABSTRACT: The first known example of a complete Riemannian manifold whose isoperimetric profile is discontinuous is given.
\par\medskip
RESUM\'E : On construit le premier exemple connu d' une vari\'et\'e riemannienne compl\`ete dont le profil isop\'erim\'etrique est discontinu.}
\end{quote}

\footnotetext[1]{S. Nardulli, Instituto de Matem\'atica, Universidade Federal de Rio de Janeiro} 
\footnotetext[2]{P. Pansu, Univ Paris-Sud, Laboratoire de
Math\'ematiques d'Orsay, Orsay, F-91405} 
\footnotetext[3]{\hskip42pt CNRS, Orsay, F-91405.}

\section{Introduction}

\subsection{The problem}

Let $M$ be a Riemannian manifold. We are concerned with the continuity of the isoperimetric profile of $M$. Given $0<v<\vol(M)$, consider all domains, i.e. smooth compact codimension $0$ submanifolds in $M$ with volume $v$. Define $I_{M}(v)$ as the least upper bound of the boundary areas of such domains. In this way, one gets a function $I_{M}:(0,\vol(M))\to\R_{+}$ called the \emph{isoperimetric profile} of $M$. 

\begin{ques}
When is the isoperimetric profile a continuous function ? 
\end{ques}

The answer is yes when $M$ is compact, Lemme 6.2 of \cite{Gallot}. S. Gallot's proof uses techniques of metric geometry. In the compact case, alternative proofs, based on the direct method of the calculus of variations, can be found in books like \cite{AmbrosioFuscoPallara}, \cite{MorganGMT}, \cite{Maggi}. 

The latter argument has been extended to the case of complete manifolds with $C^{2,\alpha}$-bounded geometry, see Theorem $1$ of \cite{NardulliFlores} and Theorem 2.2 of \cite{NardulliFlores3}. 

If one assumes existence of isoperimetric regions of every volume, one can weaken bounded geometry assumptions. It suffices to assume a lower bound on Ricci curvature and on the volumes of balls of radius 1, see Theorem 4.1 \cite{NardulliFlores}. In our opinion, it remains an open question whether the noncollapsing assumption (lower bound on the volumes of balls) can be removed or not, see Question \ref{Ques:3} below.

The isoperimetric profile is continuous also when the volume of $M$ is finite; a proof of this fact can be found in Corollary 2.4 of \cite{NardulliRusso}. 

When the ambient manifold is a non-compact homogeneous space, Hsiang showed that its isoperimetric profile is a non-decreasing and absolutely continuous function [\cite{Hsiang}, Lemma 3, Thm. 6]. 

In a recent paper, \cite{RitoreContinuity}, Manuel Ritor\'e showed that a complete Riemannian manifold possessing a strictly convex Lipschitz continuous exhaustion function has continuous and nondecreasing isoperimetric profile. Hadamard manifolds and complete non-compact manifolds with strictly positive sectional curvature belong to this class. This shows that earlier attempts to construct counterexamples using pieces of increasing negative curvature are doomed to fail. An example of a manifold with density with discontinuous isoperimetric profile has been described by Adams, Morgan and Nardulli in Prop. 2 of \cite{MorganBlog}. 

For more informations about the literature on the continuity of the isoperimetric profile the reader should consult the introduction of \cite{RitoreContinuity} and references therein. 

\subsection{The result}

\begin{thm}
\label{main}
There exists a connected non compact 3-dimensional Riemannian manifold $M$ such that $I_M$ is a discontinuous function.
\end{thm}

The proof is a modification of the treatment of Riemannian manifolds with density by Adams, Morgan and Nardulli, an account of which can be found in Frank Morgan's blog, \cite{MorganBlog}.

Start with a disjoint union of Riemannian manifolds $N=\coprod_n M_n$ such that $\vol(M_n)=1+\tau_n$ where $\tau_n>0$ tends to 0. Then $I_N(1+\tau_n)=0$. Assume that, for all $n$, $I_{M_n}(1)=I_{M_n}(\tau_n)\geq 1$. Then it is not too hard to show that $I_N(1)\geq 1$. Connecting $M_{n}$ to $M_{n+1}$ with a very thin tube produces a connected Riemannian manifold $M$ for which $I_N(1+\eps_n)$ tends to 0. Again, it is not too hard to show that $I_M(1)>0$. Therefore $I_M$ is discontinuous. 

Thus the key input is the sequence of Riemannian manifolds $M_n$ with $\vol(M_n)$ bounded and $I_{M_n}(\tau_n)$ bounded below. Adams, Morgan and Nardulli indulged themselves in introducing densities. They took for $M_n$ a tiny round sphere with a high constant density. Since volumes and boundary areas rescale differently, one can achieve $I_{M_n}(\eps_n)\geq 1$. Instead, we use nilmanifolds equipped with metrics which converge (up to rescaling) to a single Carnot-Carath\'eodory metric. The Carnot-Carath\'eodory isoperimetric profile established in \cite{Pansu83} gives a uniform lower bound for the isoperimetric profiles of such metrics.

A similar construction certainly works in any dimension $\geq 3$.

\begin{ques}\label{Ques:2}
Does there exist a 2-dimensional Riemannian manifold whose isoperimetric profile is discontinuous?
\end{ques}
\begin{ques} \label{Ques:3}
Does a manifold with Ricci curvature bounded below and admitting isoperimetric regions in every volume, have a continuous isoperimetric profile?   
\end{ques}

\section{Isoperimetry in nilmanifolds}

\subsection{Isoperimetry in the Heisenberg group}

The Heisenberg group $\H$ is the group of real upper triangular unipotent $3\times 3$ matrices,
\begin{eqnarray*}
\H=\{\begin{pmatrix}
1&x&z\\0&1&y\\0&0&1
\end{pmatrix}\,;\,x,\,y,\,z\in\R\}.
\end{eqnarray*}
Putting integer entries produces the discrete subgroup $\H_{\Z}\subset\H$. Let $dx$, $dy$, $\theta=dz-xdy$ be a basis of left-invariant forms. Let
\begin{eqnarray*}
g_{\eps}=dx^2+dy^2+\frac{1}{\eps^2}\theta^2.
\end{eqnarray*}
This is a left-invariant Riemannian metric on $\H$. As $\eps$ tends to 0, the distance $d_{\eps}$ associated to $g_\eps$ converges to the \emph{Carnot-Carath\'eodory distance}
\begin{eqnarray*}
d_c(p,q)=\inf\{\lng(\gamma)\,;\,\gamma(0)=p,\,\gamma(1)=q,\,\gamma^{*}\theta=0\}.
\end{eqnarray*}
\begin{itemize}
  \item $d_c$ has Hausdorff dimension 4, with spherical 4-dimensional measure proportional to Haar measure $\mathcal{S}^4=dxdydz$. 
    \item Smooth surfaces $S$ in $\H$ have Hausdorff dimension 3, with spherical 3-dimensional measure proportional to the measure denoted by $\mathcal{S}^3$ to be described soon (the proportionality constants are universal and will be ignored in the sequel). 
      \item Smooth curves which are transverse to the contact structure ker$(\theta)$ have Hausdorff dimension 2, with spherical 2-dimensional measure $\mathcal{S}^2$ given by integration of $\theta$ (up to sign).   \item Smooth curves tangent to ker$(\theta)$ have Hausdorff dimension 1, with spherical 1-dimensional measure $\mathcal{S}^1$ being length. 
\end{itemize}
$\mathcal{S}^3$ is locally the product $\mathcal{S}^1 \otimes \mathcal{S}^2$. Specifically, let $d\ell$ denote a (locally defined away from points where $S$ is tangent to ker$(\theta)$) unit 1-form on $S$ whose kernel is orthogonal to the trace of ker$(\theta)$ on the tangent plane to $S$. Then, up to a sign, $\mathcal{S}^3$ is obtained by integrating $d\ell\wedge\theta$. 

The Heisenberg isoperimetric inequality (\cite{Pansu83}) states that for all smooth domains $\Omega\subset \H$,
\begin{eqnarray}\label{HI}
\mathcal{S}^3(\partial \Omega)\geq \mathcal{S}^4(\Omega)^{3/4},
\end{eqnarray}
up to a universal constant that we ignore again.

Here is an alternate description of $\mathcal{S}^3$. Let $d\,\area_\eps$ denote the area induced by Riemannian metric $g_\eps$. The 1-form $\theta$ restricts to a 1-form on $S$, we denote its $g_\eps$-norm by $|\theta_{|S}|_\eps$. Then $\mathcal{S}^3$ has density $|\theta_{|S}|_\eps$ with respect to $g_\eps$-area,
\begin{eqnarray*}
d\mathcal{S}^3=|\theta_{|S}|_\eps d\,\area_\eps.
\end{eqnarray*}

Since $|\theta|_\eps=\eps$, $|\theta_{|S}|_\eps\leq\eps$, therefore $\mathcal{S}^3 (S)\leq\eps\,\area_\eps(S)$. On the other hand, the Riemannian volume element of $g_\eps$ is $\frac{1}{\eps}dxdydz$. This shows that the Heisenberg isoperimetric inequality (\ref{HI}) implies a lower bound on the isoperimetric profile of $(\H,g_\eps)$ for all $\eps>0$,
\begin{eqnarray}\label{HepsI}
I_{(\H,g_\eps)}(v)\geq \frac{1}{\eps^{1/4}}v^{3/4}.
\end{eqnarray}
This is asymptotically sharp for large volumes, but not for small volumes, where the correct asymptotics is $v^{2/3}$. Never mind, it is the dependance on $\eps$ which is most important here.

We shall not directly use inequality (\ref{HepsI}). Instead, we shall rely on inequality (\ref{HI}) to study the Carnot-Carath\'eodory isoperimetric profile of a quotient of $\H$. Only at the very end shall we return to Riemannian geometry.

\subsection{Nilmanifolds}

$\H$ possesses group automorphisms $\delta_t(x,y,z)=(tx,ty,tz)$. Let $\Gamma_t=\delta_t(\H_\Z)$ and $N_t=\Gamma_t\setminus\H$ be the quotient manifold. It inherits quotient metrics $g_\eps$, yielding Riemannian nilmanifolds $N_{t,\eps}$ of total volume equal to $\frac{t^4}{\eps}$. But is also inherits a Carnot-Carath\'eodory metric that depends only on $t$. Our first goal is to show that the Carnot-Carath\'eodory isoperimetric profile of $N_{t}$ satisfies an inequality similar to (\ref{HI}). Note that $\delta_t$ induces a homothetic map of $N_1$ onto $N_t$, so it suffices to work with one single compact space $N_1$. The volume of $N_1$ is $\mathcal{S}^{4}(N_1)=1$.

\begin{thm}\label{HIN}
There exists a constant $c$ such that the Carnot-Carath\'eodory isoperimetric profile of $N_1$ satisfies $I_{(N_1,d_c)}(v)\geq c\,\min\{v,1-v\}^{3/4}$. In other words, if $\Omega\subset N_1$ is a smooth domain of volume less that $1/2$,
\begin{eqnarray*}
\mathcal{S}^3(\partial\Omega)\geq c\,\mathcal{S}^4(\Omega)^{3/4}.
\end{eqnarray*}
\end{thm}

The method consists in cutting domains of $N_{1}$ into pieces that lift to covering spaces. Ultimately, pieces lift to $\H$ where one can apply $(\ref{HI})$. This covers cases where volume is smaller than some universal constant $v_0$. To treat domains with volume $\geq v_0$, we apply a compactness result due to \cite{Franchi-Serapioni-Serra-Cassano} (see also \cite{Leonardi-Rigot}).

\subsection{Reduction to pillars}

A first step is to cut domains into pieces called \emph{pillars} that lift to a $\Z\oplus\Z$ covering space $Z$ of $N_1$.

\begin{defi}
Let $\zeta$ denote the center of $\H_\Z$. Let us call \emph{pillar} a subset of $Z=\zeta\setminus\H$ whose projection to $\H/[\H,\H]=\R^2$ is contained in a unit square. Denote by $PI_{Z}$ the \emph{pillar profile} of $Z$, i.e.
\begin{eqnarray*}
PI_{Z}(v)=\inf\{\mathcal{S}^3(\partial\PP)\,;\,\PP \textrm{ a pillar}, \,\mathcal{S}^4(\PP)=v\}.
\end{eqnarray*} 
\end{defi}

\begin{prop}[Reduction to pillars]\label{cut}
The pillar profile of $Z$ bounds the profile of $N_{1}$ from below, with an error term,
\begin{eqnarray*}
I_{(N_{1},d_c)}(v)\geq PI_{Z}(v)-4v.
\end{eqnarray*}
\end{prop}

\proof
The coordinate functions $x$ and $y$ on $\H$ pass to a quotient $N_{1}\to \Z\setminus\R$. For $u=(s,s')\in (\Z\setminus\R)^2$, let 
\begin{eqnarray*}
G_u=\{p\in N_{1}\,;\,x(p)=s\textrm{ or }y(p)=s'\}.
\end{eqnarray*}
This is the union of two surfaces, each of which is a level set of one of the functions $x$ or $y$. The complement of $G_u$ has a cyclic fundamental group that maps isomorphically onto $\zeta$.  

Let $\Omega$ be a domain in $N_{1}$. By the coarea formula,
\begin{eqnarray*}
\mathcal{S}^4(\Omega)=\int_{\Z\setminus\R}\mathcal{S}^3(x^{-1}(s)\cap\Omega)\,ds.
\end{eqnarray*}
This coarea formula follows from the fact that the volume element (viewed as a 4-form) splits,
\begin{eqnarray*}
d\,\mathcal{S}^4=dx\wedge dy\wedge\theta=dx\wedge d\,\mathcal{S}^3 ,
\end{eqnarray*}
since $dy\wedge\theta=d\,\mathcal{S}^3$ along the fibers of $x$ (one can take $d\ell=dy$ globally). 

The same inequality holds with $x$ replaced with $y$. This shows that there exists $u=(s,s')\in (\Z\setminus\R)^2$ such that 
\begin{eqnarray*}
\mathcal{S}^3(x^{-1}(s)\cap\Omega)\leq\mathcal{S}^4(\Omega), \quad \mathcal{S}^3(y^{-1}(s')\cap\Omega)\leq\mathcal{S}^4(\Omega),
\end{eqnarray*}
and thus
\begin{eqnarray*}
\mathcal{S}^3(G_u\cap\Omega)\leq 2\mathcal{S}^4(\Omega).
\end{eqnarray*}
The complement $\Omega\setminus G_u$ lifts to the cyclic covering space $Z$. Pick some lift. Its closure $\PP$ is a pillar. Indeed, on $\PP$, the real valued functions $x$ and $y$ take values in intervals of length $1$. The boundary of $\PP$ consists of a part that isometrically and injectively maps to $\partial\Omega$, and of a part that maps 2-1 to $G_u\cap\Omega$. Therefore
\begin{eqnarray*}
\mathcal{S}^3(\partial\PP)\leq\mathcal{S}^3(\partial\Omega)+2\mathcal{S}^3(G_u\cap\Omega)\leq\mathcal{S}^3(\partial\Omega)+4\mathcal{S}^4(\Omega).
\end{eqnarray*}
If $\mathcal{S}^4(\Omega)=v$, this shows that
\begin{eqnarray*}
I_{(N_{1},d_c)}(v)\geq PI_{Z}(v)-4v.\qed
\end{eqnarray*}

\subsection{Treatment of pillars}

\begin{prop}[Treatment of pillars]\label{recut}
The profile of $\H$ bounds the pillar profile of $Z$ from below, with an error term,
\begin{eqnarray*}
PI_{Z}(v)\geq I_{\H}(v)-8v.
\end{eqnarray*}
\end{prop}

\proof
Let $\PP\subset Z$ be a pillar. We can assume that its projection to $\R^2$ is contained in $\{0\leq x\leq 1\}$. Its inverse image $\tilde{\PP}$ in $\H$ is a $\zeta$-invariant subset with small projection in $\R^2$. Again, we cut $\tilde{\PP}$ into logs of height $1$ using level sets of the $z$ function. This time, we split the volume element as
\begin{eqnarray*}
d\mathcal{S}^4=dx\wedge dy\wedge dz=dz\wedge(dx\wedge dy)=dz\wedge\frac{1}{|x|}d\,\mathcal{S}^3\geq dz\wedge d\,\mathcal{S}^3.
\end{eqnarray*}
We have used the expression $d\,\mathcal{S}^3 =|x|\,dx\,dy$ for the measure induced on horizontal planes $\{z=s\}$. On such surfaces, one can take $d\ell=dx$ globally, whence $d\,\mathcal{S}^3 =\pm dx\wedge\theta=|x|\,dx\,dy$. The coarea formula gives
\begin{eqnarray*}
\mathcal{S}^4(\PP)&=&\mathcal{S}^4(\tilde{\PP}\cap\{0\leq z\leq 1\})\\
&=&\int_{0}^{1}\left(\int_{\tilde{\PP}\cap\{z=s\}}\frac{1}{|x|}d\,\mathcal{S}^3\right)\,ds\\
&\geq&\int_{0}^{1}\mathcal{S}^3(\tilde{\PP}\cap\{z=s\})\,ds.
\end{eqnarray*}
There exists $s\in[0,1]$ such that
\begin{eqnarray*}
\mathcal{S}^3(\tilde{\PP}\cap\{z=s\})\leq\mathcal{S}^4(\PP).
\end{eqnarray*}
Set $\Omega'=\tilde{\PP}\cap\{s\leq z\leq s+1\}$. Then 
\begin{eqnarray*}
\mathcal{S}^3(\partial\Omega')\leq \mathcal{S}^3(\partial\PP)+2\mathcal{S}^4(\PP).
\end{eqnarray*}
If $\PP$ has volume $v$, this leads to
\begin{eqnarray*}
PI_{Z}(v)\geq I_{\H}(v)-2v.
\end{eqnarray*}

\subsection{Profile of $(N_1,d_c)$}

\begin{prop}[Carnot-Carath\'eodory isoperimetric inequality for small volumes]\label{CC}
If $v\leq 12^{-4}$,
\begin{eqnarray*}
I_{(N_{1},d_c)}(v)\geq \frac{1}{2}v^{3/4}.
\end{eqnarray*}
\end{prop}

\proof
Combined with Propositions \ref{cut} and \ref{recut}, the Heisenberg isoperimetric inequality (\ref{HI}) yields
\begin{eqnarray*}
I_{(N_{1},d_c)}(v)\geq v^{3/4}-4v-2v=v^{3/4}(1-6\, v^{1/4})\geq\frac{1}{2}v^{3/4},
\end{eqnarray*}
since $v\leq 12^{-4}$. \qed

\subsection{Proof of Theorem \ref{HIN}}

There is a notion of Carnot-Carath\'eodory perimeter, an appropriate topology for which $\mathcal{S}^4$ is continuous and the perimeter (which coincides with $\mathcal{S}^3$ for smooth domains) lower semi-continuous, and a compactness theorem for sets of bounded perimeter in a compact Carnot manifold, \cite{Franchi-Serapioni-Serra-Cassano}. This implies that the Carnot-Carath\'eodory isoperimetric profile $I_{(N_1,d_c)}$ is positive on $(0,1)$ and lower semi-continuous. Therefore, there exists $\eta>0$ such that $I_{(N_1,d_c)}\geq \eta$ on $[12^{-4},1-12^{-4}]$. Set $c=\min\{\frac{1}{2},2^{3/4} \eta\}$. Then $I_{(N_1,d_c)}(v)\geq\eta=c(\frac{1}{2})^{3/4}\geq c\,v^{3/4}$ for every $v\in[12^{-4},\frac{1}{2}]$. On the other hand, Proposition \ref{CC} shows that $I_{(N_1,d_c)}(v)\geq c\,v^{3/4}$ for all $v\in[0,12^{-4}]$. \qed

Note that the proof does not provide an effective constant $c$.

\subsection{Riemannian profile}

\begin{cor}
Let $N_{t,\eps}$ denote the quotient $(\delta_t(\H_\Z))\setminus\H$ equipped with the Riemannian metric induced by $g_\eps$. The isoperimetric profile of $N_{t,\eps}$ satisfies
\begin{eqnarray*}
I_{N_{t,\eps}}(v)\geq \frac{c}{\eps^{1/4}}\,\min\{v,\frac{t^4}{\eps}-v\}^{3/4}.
\end{eqnarray*}
\end{cor}

\proof
The homothetic map $N_1\to N_t$ induced by the automorphism $\delta_t$ transports the inequality of Theorem \ref{HIN} to $N_t$ without any change but the fact that $\mathcal{S}^4(N_t)=t^4$ replaces $1$. The Riemannian volume element of $N_{t,\eps}$ is $\frac{1}{\eps}\mathcal{S}^4$, the Riemannian area induced on surfaces satisfies $\eps\,\area \geq \mathcal{S}^3$. This leads to the indicated dependence on $\eps$ in the isoperimetric profile of $N_{t,\eps}$. \qed

\section{Proof of Theorem \ref{main}}

\subsection{The case of a disjoint union of nilmanifolds}

\begin{prop}\label{disjoint}
Let $\tau_n=\frac{1}{n}$, $\eps_n=\tau_n^3$ and $t_n=\tau_n^{3/4}(1+\tau_n)^{1/4}$. Let $N=\coprod_{n}N_{t_n,\eps_n}$. Then, for all $v\in[\frac{1}{16},1]$, $I_N(v)\geq \frac{c}{8}$, where $c$ is the constant of Theorem \ref{HIN}.
\end{prop}

\proof
By construction, $\vol(N_{t_n,\eps_n})=1+\tau_n$.
Let $\Omega$ be a domain in $N$ with $\vol(\Omega)=v$. Write $\Omega=\coprod_n \Omega_n$ where $\Omega_n\subset N_{t_n,\eps_n}$ has volume $v_n$, $\sum_{n=1}^{\infty}v_n=v$. 

If some $v_n$ satisfies $v_n\geq \frac{1}{2}(1+\tau_n)$, then 
\begin{eqnarray*}
\area(\partial\Omega_n)&\geq&\frac{c}{\eps_n^{1/4}}(1+\tau_n-v_n)^{3/4}\\
&\geq&\frac{c}{\eps_n^{1/4}}\tau_n^{3/4}= c,
\end{eqnarray*}
so 
\begin{equation}\label{Eq:disjoint}
\area(\partial\Omega)\geq c
\end{equation} in this case. 

Otherwise, for all $n\geq 1$,
\begin{eqnarray*}
\area(\partial\Omega_n)\geq\frac{c}{\eps_n^{1/4}}v_n^{3/4}\geq c\, v_n^{3/4}.
\end{eqnarray*}
We use the concavity inequality
\begin{eqnarray*}
a^{\alpha}+b^{\alpha}\geq (a+b)^{\alpha},
\end{eqnarray*}
valid for all $0\leq\alpha\leq 1$, $a\geq 0$ and $b\geq 0$. This gives
\begin{eqnarray*}
\area(\partial\Omega)&=&\sum_{n=1}^{\infty}\area(\partial\Omega_n)\\
&\geq&c\sum_{n=1}^{\infty}v_n^{3/4}\\
&\geq&c\left(\sum_{n=1}^{\infty}v_n\right)^{3/4}=(\frac{1}{16})^{3/4}c=\frac{c}{8}.\qed
\end{eqnarray*}

\subsection{Connecting manifolds}
\proof
We construct a noncompact manifold that has the shape of an infinite pearl necklace, adjusting suitable parameters carefully.
Let $0<\tau_n<1$ be the sequence of positive real numbers chosen in the proof of Proposition \ref{disjoint}. Pick another sequence of volumes $w_n<1$, such that
\begin{equation}\label{Eq:Main1}
\sum_n w_n<\frac{1}{2},
\end{equation}
and a sequence of areas $a_n>0$ such that
\begin{equation}\label{Eq:Main}
\sum_n a_n<\frac{c}{16},
\end{equation}
where $c$ is the constant of Theorem \ref{HIN}.

The manifolds $N_{t_n,\eps_n}$ that we want to connect to obtain our counterexample $M$, are like in Proposition \ref{disjoint}, in particular we retain here that $V(N_{t_n,\eps_n})=1+\tau_n$. Take two small disjoint balls $B_{n,1}, B_{n,2}$ inside $N_{t_n,\eps_n}$ whose boundaries have total area $\leq a_n$. Arrange that $B_{n,2}$ and $B_{n+1,1}$ be nearly isometric. Put $\tilde{N}_n:=N_{t_n,\eps_n}\setminus\left(B_{n,1}\mathring{\cup}B_{n,2}\right)$. 

Consider tubes or cylinders $T_n$ of the form $T_n:=(S^2(1)\times [0, 1], g_n)$, where the metrics $g_n$ are chosen in such a way that $V(g_n)\leq w_n$ and they glue together into a smooth metric on the connected sum $M_n:=\tilde{N}_n\#T_n$ where the gluing is done along $i_n(S^2(1)\times\{0\})\cong\partial B_{n,2}$. Now consider 
\begin{equation}
(M,g):=M_1\#M_2\#\cdots\#M_{n}\#M_{n+1}\#\cdots
\end{equation} 
where $M_n$ and $M_{n+1}$ are glued together along the boundaries $i_n(S^2(1)\times\{1\})\cong\partial B_{(n+1),1}$, where $i_n:T_n\rightarrow M$ is the isometric embedding associated to our construction. 

We show that the right limit $I_M(1+)$ vanishes. Consider domains $D_n:=\tilde{N}_n$, we get $V(D_n)=1+\tau_n-\tilde{v}'_{n-1}-\tilde{v}'_n=1+\alpha_n$, with $\alpha_n\rightarrow 0$, $\varepsilon'_n:=A(\partial D_n)=A_g(\partial B_{n,2}\mathring{\cup}\partial B_{n+1,1})\rightarrow 0$. This implies readily \begin{equation}0\leq\lim_{n\rightarrow+\infty} I_M(1+\alpha_n)\leq\lim_{n\rightarrow+\infty} A(\partial D_n)=0.\end{equation}

We show that $I_M(1)>0$. Let $\Omega$ be a domain in $M$ such that $V(\Omega)=1$. Write $\tilde{\Omega}:=\mathring{\bigcup}\tilde{\Omega}_n$, where $\tilde{\Omega}_n:=\Omega\cap\tilde{N}_n$. Then 
$$V(\tilde{\Omega})\geq 1-\sum_n w_n \geq \frac{1}{2}.$$
According to Proposition \ref{disjoint},
$$A(\partial\tilde{\Omega})\geq \frac{c}{8}.$$
Since, for all $n$,
 \begin{eqnarray*}
\partial\tilde{\Omega}_n=((\partial\Omega)\cap\tilde{N}_n) \mathring{\cup} (\Omega\cap \partial\tilde{N}_n),
\end{eqnarray*}
\begin{eqnarray*}
A(\partial\tilde{\Omega}_n)-A((\partial\Omega)\cap\tilde{N}_n)\leq A_g(\partial B_{n,2}\mathring{\cup}\partial B_{n,1})\leq a_n,\end{eqnarray*}
thus
\begin{eqnarray*}
A(\partial\Omega)\geq A(\partial\tilde{\Omega})-\sum_{n}a_n\geq \frac{c}{8}-\frac{c}{16}=\frac{c}{16}.
\end{eqnarray*}
This show that $I_M(1)\geq\frac{c}{16}$.
 
 This concludes the proof of Theorem \ref{main}.\qed
 
\newpage
      \markboth{References}{References}
      \bibliographystyle{alpha}
      \bibliography{Discontinuous}
      \addcontentsline{toc}{section}{\numberline{}References}

Keywords : Isoperimetric inequality, Nilmanifold, Carnot-Carath\'eodory metric. 

Mathematics Subject Classification : 
53C20, 
49Q20, 

\vskip1cm
\noindent
Stefano Nardulli\\
Instituto de Matem\'atica\\
Universidad Federal de Rio de Janeiro (Brazil)\\
\smallskip\noindent
{\tt nardulli@im.ufrj.br}\\
http://www.im.ufrj.br/nardulli
\par\medskip\noindent
Pierre Pansu\\
Laboratoire de Math\'ematiques d'Orsay, UMR 8628 du C.N.R.S.\\
B\^atiment 425, Universit\'e Paris-Sud - 91405 Orsay (France)\\
\smallskip\noindent
{\tt\small Pierre.Pansu@math.u-psud.fr}\\
http://www.math.u-psud.fr/$\sim$pansu

\end{document}